\documentclass[10pt, reqno]{amsart}
\usepackage[arrow,matrix,curve]{xy}

\newtheorem{theorem}{Theorem}[section]
\newtheorem{proposition}[theorem]{Proposition}

\newtheorem{corollary}[theorem]{Corollary}

\theoremstyle{definition}

\newtheorem{question}[theorem]{Question}

\numberwithin{equation}{section} 
\numberwithin{figure}{section}

\begin{document}

\author[C.~Croke]{Christopher B. Croke$^+$} \address{ Department ofbb:
Mathematics, University of Pennsylvania, Philadelphia, PA 19104-6395
USA} \email{ccroke@math.upenn.edu} \thanks{$^+$Supported by NSF grant
DMS 02-02536 and the Max-Planck Inst. Bonn.}


\title[hemispheres] {A synthetic characterization of the hemisphere}


\begin{abstract}
We show that round hemispheres are the only compact 2 dimensional Riemannian manifolds (with or without boundary) such that almost every pair of complete geodesics intersect once and only once.  We prove this by establishing a sharp isoperimetric inequality for surfaces with boundary such that every pair of geodesics have at most one interior intersection point.
\end{abstract}

\maketitle

\section{Introduction}
In this short note we prove a sharp isoperimetric inequality for 2 dimensional Riemannian manifolds $(M,\partial M,g)$ with boundary $\partial M$.  
The result will apply to surfaces such that every pair of complete geodesics have at most one interior intersection point.  So in particular it applies to surfaces where all geodesics are minimizing.  By a ``complete geodesic'' we mean a curve satisfying the geodesic equation which either runs from a boundary point to a boundary point or is infinitely long and the boundary points of the curve are precisely the intersection points of the curve with the boundary.  They can a-priori have 0, 1, or 2 boundary points.  However the intersection condition (along with the finiteness of the area) will quickly allow us to conclude that all geodesics are compact segments between two boundary points.  We should point out that we are ruling out closed geodesics since they are interpreted to intersect any geodesic infinitely often if they intersect at all.  (Note that we would have to include round projective planes if one interpreted such an intersection as a single intersection.)

\begin{theorem}
\label{isoperimetric}
Let $(M,\partial M,g)$ be a two dimensional, finite area, Riemannian manifold with (possibly empty) boundary, $\partial M$, such that every pair of complete geodesics have at most one interior point of intersection.  Then
$$L(\partial M)^2\geq 2\pi Area(M).$$
Further equality holds if and only if every pair of geodesics intersect.
\end{theorem}

This is a new proof of (and an extension of) the inequality in the following case of the isoperimetric inequality Theorem 11 of \cite{Cr}:  

\begin{theorem}
\label{isoperimetric2}
Let $(M,\partial M,g)$ be a 2 dimensional Riemannian manifold such that all geodesics minimize until they hit the boundary.  Then
$$L(\partial M)^2\geq 2\pi Area(M).$$
Further equality holds if and only if $M$ is isometric to a round hemisphere.
\end{theorem}

In fact Theorem 11 of \cite{Cr} says much more than this, and applies in higher dimensions.  

The interest in Theorem \ref{isoperimetric} is twofold.  First, the equality cases in the above give us the characterization of hemispheres: 

\begin{corollary}
\label{characterization}
Round hemispheres are the only finite area, 2 dimensional, Riemannian manifolds (with or without boundary) such that almost every pair of geodesics intersect once and only once. \end{corollary}

In the infinite volume case both Theorem \ref{isoperimetric} and Corollary \ref{characterization} are false because we would have to include the flat plane.  In fact, below is a (probably quite difficult) question that has been open for some time. 

\begin{question}
\label{euclid}
Is the Euclidean metric the only complete metric on the plane with no conjugate points such that for every line $l$ (complete geodesic) and point $p$ not on $l$ there is a unique line through $p$ parallel to (i.e. not intersecting) $l$?
\end{question}

The answer would be yes if there was a positive answer to:

\begin{question}
\label{synth}
Is the Euclidean metric the only complete metric on the plane such that almost every pair of geodesics intersect in exactly one point?
\end{question}

Some results on Question \ref{euclid} can be found in \cite{Bu-Kn}.  Question \ref{synth} was also posed by Knieper in an MSRI problem session in the early nineteen-nineties.

We should point out that although there can be no direct higher dimensional analogue of Corollary \ref{characterization} the closest would be:  {\em n-dimensional hemispheres are the only compact manifolds with boundary such that all complete geodesics have the same length.}  This was proved in \cite{Ba} using the Blaschke conjecture for spheres. It also follows from Santal\'o's formula (see below) and the equality case of Theorem 11 of \cite{Cr}.  

The second reason for interest in Theorem \ref{isoperimetric} is that the proof tells us that the isoperimetric deficit, $L(\partial M)^2 - 2\pi Area(M)$ is $\frac 1 4$ the measure of the set of pairs of geodesics that do not intersect (i.e. $\frac 1 4 Vol(\Gamma\times \Gamma - I)$).  This will follow from the Proposition \ref{area} of the next section.  In particular, if $M$ in addition has nonpositive curvature then (see \cite{We}) it satisfies the classical isoperimetric inequality, $L(\partial M)^2\geq 4\pi Area(M)$, and hence for such metrics $\frac 1 4 Vol(\Gamma\times \Gamma - I) \geq 2\pi Area(M)$ with equality holding only for the flat disc.

\section{Proofs}

We will assume throughout this section that $M$ is a surface of finite area with (possibly empty) boundary $\partial M$.  We can assume that the boundary is compact for if not it would have infinite length and the inequality would be clearly true (and not sharp). The unit tangent bundle, $UM$, of $M$ has a natural measure, $du$, that is locally a product measure.  In particular $Vol(UM)=2\pi Area(M)$.  Let $\tau$ be a curve in $M$ with arclength $s$ and unit normal $\nu$.  We can let $\{(\theta,s)\}$ parameterize the set of geodesics $\gamma$ that intersect $\tau$, by $\gamma(0)= \tau(s)$ and the unit tangent $\gamma'(0)$ makes angle $\theta$ with $\nu$.  For all $t$ such that $\gamma(t)$ is defined we let $\{(\theta,s,t)\}$ correspond to the unit vector $u=\gamma'(t)$.  Then Santal\'o's formula tells us that $du=|cos(\theta)|d\theta ds dt.$  We will be using this in a number of ways.

First consider metrics where any two geodesics intersect at most once in the interior.  Note that this will follow (by continuity) if we assume that almost every pair intersects at most once.  Some pairs may intersect more than once but then only at boundary points (as is the case for round hemispheres).  We claim that all complete geodesics are segments between boundary points.  If not then let our infinite length geodesic play the role of $\tau$ above.  We first note that $\tau:[0,\infty)\to M$ does not have a limit point on the boundary since if so a geodesic leaving normal to the boundary at that point would intersect $\tau$ more than once.  In particular there is an $\epsilon>0$ such that $d(\tau[1,\infty),\partial M)>\epsilon$.
Since every geodesic $\gamma$ that intersects $\tau$ intersects it only once the parametrization $(\theta,s,t)$ for $1\leq s<\infty$ and $-\epsilon<t<\epsilon$ of $UM$ is one-to-one (though not onto).  Hence $Vol(UM)\geq 2\pi\times \infty\times 2\epsilon$ a contradiction.  

We now consider the case where all geodesics go between boundary points (without assuming geodesics intersect at most once).  It is easy to describe the (standard) measure space, $\Gamma$, of complete unit speed oriented geodesics by using the boundary $\partial M$ as our $\tau$ above.  The measure is $d\gamma=|cos(\theta)|d\theta ds$.  Note that the parametrization $(\theta,s,t)$ of $UM$ is one-to-one and onto if $-\frac \pi 2\leq \theta \leq \frac \pi 2$, $s$ is an arc length parametrization of the boundary (which is not assumed to be connected) and $0\leq t\leq L(\gamma)$ (where $\gamma$ is the geodesic determined by $(\theta,s)$).  
In particular, $2 \pi Area(M)=Vol(UM)=\int_\Gamma L(\gamma)d\gamma$.
   
Now consider a curve $\tau$.  We claim that $\int_\Gamma i(\tau,\gamma)d(\gamma)=4L(\tau)$ where $i(\tau,\gamma)$ represents the number of times the geodesic $\gamma$ intersects $\tau$.  (This is known as Crofton's formula.)  It follows directly from Santal\'o's formula when we note that the parametrization $(\theta,s)$ of $\Gamma$ counts each $\gamma$ exactly as often as $i(\tau,\gamma)$ and hence $4L(\tau)=\int\int|cos(\theta)|d\theta ds=\int_\Gamma i(\tau,\gamma)d(\gamma)$.
In particular (since each geodesic hits $\partial M$ twice) we see that $Vol(\Gamma)=2L(\partial M)$.

\begin{proposition}
\label{area}
Let $(M,\partial M,g)$ be such that every complete geodesic hits the boundary at both ends.  Then
$$8\pi Area(M)=\int_{\Gamma \times \Gamma}i(\gamma_1,\gamma_2)d\gamma_1 d\gamma_2.$$

In particular, if every pair of geodesics intersect at at most one interior point then 
$8\pi Area(M)=Vol(I)$ where $I\subset \Gamma \times \Gamma$ represents the subset of geodesic pairs that intersect. 
\end{proposition}
 
\proof:
$$\int_{\Gamma \times \Gamma}i(\gamma_1,\gamma_2)d\gamma_1 d\gamma_2=
\int_\Gamma\Big \{\int_\Gamma i(\gamma_1,\gamma_2)d\gamma_1\Big \}d\gamma_2=
$$$$=\int_\Gamma 4L(\gamma_2)d\gamma_2=8\pi Area(M).
$$
\qed

Note that Theorem \ref{isoperimetric} follows immediately from this since $Vol(I)\leq Vol(\Gamma \times \Gamma)= 4L(\partial M)^2$ and equality holds if and only if $\Gamma \times \Gamma - I$ has measure 0.  By continuity this implies that every pair of geodesics must intersect.

We now prove Corollary \ref{characterization}.  First note that there can be no pair of conjugate points along any geodesic segment except for the endpoints, for if so it is easy to see that there are nearby geodesics that intersect twice in the interior and hence a set of positive measure of geodesic pairs that intersect more than once.  We also note that the exponential map from any interior point is one-to-one.  Although (at least a-priori) this may not mean all geodesics minimize, this is exactly what is used in the proof of Theorem 11 of \cite{Cr}.  Thus we can apply the equality case in that theorem to get the Corollary.

\end{document}